\documentclass[12pt]{amsart}
\usepackage{latexsym}
\newenvironment{theorem*}
{\begin{trivlist}\item[\hskip%
\labelsep{\bf{Theorem.}}]\it}%
{\end{trivlist}}
\newenvironment{lemma*}
{\begin{trivlist}\item[\hskip%
\labelsep{\bf{Lemma.}}]\it}%
{\end{trivlist}}
\newenvironment{claim*}
{\begin{trivlist}\item[\hskip%
\labelsep{\itshape{Claim.}}]\it}%
{\end{trivlist}}
\newenvironment{definition*}
{\begin{trivlist}\item[\hskip%
\labelsep{\bf{Definition.}}]\rm}%
{\end{trivlist}}
\newenvironment{remark*}
{\begin{trivlist}\item[\hskip%
\labelsep{\bf{Remark.}}]\rm}%
{\end{trivlist}}
\def\OX{\mathcal{O}_n(X)}
\def\On1{\mathcal{O}_n(X_1)}
\begin{document}
\author[F.~Torres]{Fernando Torres}\thanks{The author is supported by a 
grant from the International Atomic Energy Agency and UNESCO}
\title[Constellations of Weierstrass points]{On the 
constellations of Weierstrass points}
\address{ICTP, Mathematics Section, P.O. Box 586, 34100, Trieste - Italy}
\email{feto@ictp.trieste.it}
\begin{abstract}
We prove that the constellation of Weierstrass points characterizes the 
isomorphism-class of double coverings of curves of genus large enough.
\end{abstract}
\maketitle
\noindent {\bf 1.} Let $X$ be a projective, irreducible, non-singular 
algebraic curve defined over an algebraically closed field $k$ of 
characteristic $p$. Let 
$n\ge 1$ be an integer and $C_X$ a canonical divisor of $X$. The 
pluricanonical linear system $|nC_X|$ defines a nondegenerate morphism
$$
\pi_n: X \to \mathbb P^{N(n)},
$$
where $N(1)=g-1$, and $N(n)= (2n-1)(g-1)-1$ for $n\ge 2$. To any $P\in X$ 
we then associate the sequence of multiplicities 
$$
\{v_P(\pi^*_n(H)): H\ {\rm hyperplane} \subseteq \mathbb P^{N(n)}\}=
\{\epsilon_0(P)<\epsilon_1(P)<\ldots<\epsilon_{N(n)}(P)\}.
$$
Such a sequence is the same for all but finitely many points 
(cf. \cite[III.5]{F-K}, \cite[Prop.3]{Lak}, \cite[\S1]{S-V}). These finitely 
many points
are the so called $n$-Weierstrass points of $X$. There exists a divisor 
$\mathcal{W}_n$ on $X$ whose support is the set of $n$-Weierstrass 
points and satisfies the property below. Let denote by 
$v_P(\mathcal{W}_n)$ the coefficient of 
$\mathcal{W}_n$ in $P$ (called the $n$-Weierstrass weight at $P$). Then
\begin{equation}\label{hur}
\omega_n:= {\rm deg}(\mathcal{W}_n)=\sum_{P}v_P(\mathcal{W}_n) = 
\sum_{i=0}^{N(n)}\epsilon_i (2g-2) + n(2g-2)(N(n)+1),
\end{equation}
where $\epsilon_0<\epsilon_1<\ldots<\epsilon_{N(n)}$ denotes the 
sequence at a generic point (\cite[III.5]{F-K}, \cite[Thm.6]{Lak}, 
\cite[p.6]{S-V}). One has $\epsilon_i(P)\ge \epsilon_i$ for 
each $i$ and for each $P$ (\cite[III.5]{F-K}, \cite[Prop.3]{Lak}, 
\cite[p.5]{S-V}).

Let $P_1,\ldots,P_{\omega_n}$ be the $n$-Weierstrass points (counted with 
multiplicity according to their $n$-Weierstrass weights).
\begin{definition*}
The orbit $\OX$ of $(\pi_n(P_1),\ldots,\pi_n(P_{\omega_n}))\in 
(\mathbb P^{N(n)})^{\omega_n}$ under the action of the product of 
the symmetric group 
$S_{\omega_n}$ and the projective linear group $PGL(N(n)+1)$ is called 
the constellation of $n$-Weierstrass points of $X$. $\mathcal{O}_1(X)$ 
is called the constellation of Weierstrass points of $X$. 
\end{definition*}

Let $X$, $X_1$ be curves as above of genus $g$. Then clearly $\OX =\On1$ 
if $X$ is isomorphic to $X_1$. In this note we are interested in the 
converse: 
\begin{equation*}
  \OX =\On1 \ \mbox{ for some}\ n\ge 1 \Rightarrow X\ \mbox{is 
isomorphic to\ } X_1\ ?\tag{$P$}
\end{equation*}
This problem was studied by Pflaum \cite{Pf} for $p=0$. He showed  
that $(P)$ is true in the following cases :
\begin{enumerate}
\item[(i)] If $n\ge 2$;
\item[(ii)] If $n=1$, and $2\le g\le 15$;
\item[(iii)] If $n=1$, and $X$ and $X_1$ are $\gamma$-hyperelliptic curves 
(that 
is double coverings of curves of genus $\gamma$) with $g\ge 2$ and  
$\gamma\in \{0,1,2\}$.
\end{enumerate}

For the cases $n=2$, $g=2$; $n=1$, $g=3,4$; and $n=1$, $X$, $X_1$ 
hyperelliptic curves the proof is given by a direct application of the 
definition above. For the remaining cases Pflaum used the following 
argument. By means of a lower bound on the number 
$W_n$ of $n$-Weierstrass points, he stated (Corollary 2.6, loc. cit.) a 
sufficient condition to have an affirmative answer to the problem (this 
condition holds regardless the characteristic of $k$).   
Define  the number $N(g,n)$ by $N(6,1)=25$, $N(g,1)= {\rm 
max}\{3g+6, 4g-4\}$ for $g\neq 6$, and $N(g,n)= 4n(g-1)$ for $n\ge 2$.  
Pflaum then showed that $(P)$ is true provided  
\begin{equation}\label{bound}
W_n>N(g,n).
\end{equation}
A way of getting a lower bound on $W_n$ is by bounding from above 
$v_P(\mathcal{W}_n)$, and then by using (\ref{hur}). One has 
\begin{equation}\label{weight}
v_P(\mathcal {W}_n) \ge \sum_{i=1}^{N(n)}(\epsilon_i(P)-\epsilon_i)
\end{equation}
for each $P$, and equality holds if
$$
{\rm det}\biggl(\binom{\epsilon_i(P)}{\epsilon_j}\biggr)\not\equiv 0 \pmod{p}
$$
(\cite[p.239]{Lak}, \cite[Thm.1.5]{S-V}). A 
curve $X$ is called {\it classical} (with respect to the pluricanonical 
linear system $|nC_X|$) if $\epsilon_i=i$ for each $i$. This is the case 
if $p=0$ or $p>n(2g-2)$, and here one also has the equality in 
(\ref{weight}) for each 
$P$  (\cite[III.5]{F-K}, \cite[Corollary 1.8]{S-V}, \cite[Thm.11]{Lak}).

Let $X$ and $X_1$ be classical curves of genus $g$, and suppose that we have 
equality in 
(\ref{weight}) for each point of $X$ and $X_1$. If $n=1$ and $5\le g\le 15$, 
Pflaum 
showed (\ref{bound}) by direct computations. If $n\ge2$, then Homma and 
Ommori \cite{H-O} noticed 
that $v_P(\mathcal{W}_n)\le g(g+1)/2$. Hence by means of (\ref{hur}) Pflaum  
obtained (\ref{bound}). Now if $X$ and $X_1$ are $\gamma$-hyperelliptic 
curves of genus large enough with $\gamma\in\{1,2\}$ and if $p=0$ or 
$p>2g-2$, by using results due to Garcia \cite{G} and 
Kato \cite{K} Pflaum can bound from above  $v_P(\mathcal{W}_1(P))$ and 
obtain (\ref{bound}). 

The aim of this note is to extend Pflaum's result (iii) above to 
$\gamma$-hyperelliptic curves of genus large enough ($\gamma\ge 3$) and 
whenever $p=0$ or $p>2g-2$. We will show that 
such curves satisfy (\ref{bound}) and to do that we use some results 
concerning  
Weierstrass weights in \cite{To1} and \cite{To2}. We show
\begin{theorem*}
Let $X$ and $X_1$ be $\gamma$-hyperelliptic curves curves of genus 
$g\ge 9\gamma-17+\frac{43\gamma -20}{2\gamma^2+\gamma-1}$ with 
$\gamma\ge 3$. Assume that $p=0$ or $p>2g-2$. Then 
$X$ and $X_1$ are isomorphic provided  $\mathcal{O}_1(X) = 
\mathcal{O}_1(X_1)$. 
\end{theorem*}

In general one cannot expect to fulfil condition (\ref{bound}) 
for $0<p\le 2g-2$, because in this case there exist curves with small number 
of Weierstrass points. For instance there exist curves with just one 
Weierstrass point (see \cite[\S6]{Lak}).
\medskip

\noindent {\bf 2.} Let 
$$
\pi:X\to \tilde X
$$
be a double covering of curves of genus $g$ and $\gamma\ge 3$ respectively. 
Let $P\in 
X$, and set $w(P):= v_P(\mathcal{W}_1)$. The key point of the proof is 
the fact that we can bound from above $w(P)$ by considering the following 
three cases:
\begin{enumerate}
\item[(I)] $P$ is a ramified point of $\pi$ such that $\pi(P)$ is a 
Weierstrass point of $\tilde X$.
\item[(II)] $P$ is a ramified point of $\pi$ such that $\pi(P)$ is not a 
Weierstrass point of $\tilde X$.
\item[(III)] $P$ is not a ramified point of $\pi$.
\end{enumerate}
\begin{lemma*} Let $X$ be a $\gamma$-hyperelliptic 
curve of genus $g$ ($\gamma\ge 3$). Let assume that $p=0$ or $p>2g-2$, and 
let $P\in X$. 
\begin{enumerate}
\item[(i)] If $P$ is as in (I), then
$$
w(P)\le c_1:= \binom{g-2\gamma}{2} + 2\gamma^2.
$$
\item[(ii)] If $P$ is as in (II), then
$$
w(P) \le c_2:= \binom{g-2\gamma}{2} + 4\gamma-4.
$$
\item[(iii)] If $P$ is as in (III) and $g\ge 2\gamma$, then
$$
w(P)\le c_3:= \max\{(2\gamma-1)g-(\gamma-1)(2\gamma+1), 
2\gamma g - 2\gamma(4\gamma-1)\}.
$$
\end{enumerate}
\end{lemma*}
\noindent {\bf 3. Proof of the Theorem.} Since $g\ge 
2\gamma+2$, every ramified point of $\pi$ 
is a Weierstrass point (see \S4). Let $t$ denote the number of 
points of type (I). Then by the lemma we have:
\begin{align*}
{\rm deg}(\mathcal{W}_1)=g^3-g & \le 
tc_1+ (2g-4\gamma+2-t)c_2+(W_1-2g+4\gamma-2)c_3\\
                               & = (c_1-c_2)t+ (2g-4\gamma+2)c_2+
(W_1-2g+4\gamma-2)c_3.
\end{align*}
Then by noticing that $t\le \min\{\gamma^3-\gamma,2g-4\gamma+2\}$, we find
\begin{equation}\label{w1}
(W_1-2g+4\gamma-2)c_3 \ge 6\gamma g^2 -16\gamma^2g 
+16\gamma^3-4\gamma^2-2\gamma.
\end{equation}
If $g\ge 6\gamma^2-\gamma+1$, then $c_3=2\gamma(g-4\gamma+1)$ and from 
(\ref{w1}) we get 
$$
W_1\ge 5g-1+\frac{24\gamma^2-10\gamma-4}{g-4\gamma+1}>N(g,1).
$$
Now suppose that $c_3=(2\gamma-1)g-(2\gamma^2-\gamma-1)$. From (\ref{w1}) 
we find that $W_1>N(g,1)$ if
$$
(2\gamma-1)(2\gamma+2)g - (36\gamma^3-50\gamma^2+34\gamma-6)+
\frac{24\gamma^5-40\gamma^4+
22\gamma^3+4\gamma}{(2\gamma-1)g-(2\gamma^2-\gamma-1)}>0.
$$
This is satisfied provided $g\ge 9\gamma 
-17+\frac{43\gamma-20}{2\gamma^2+\gamma-1}$. 
\medskip

\noindent {\bf 4. Proof of the Lemma.} First we recall some properties of 
Weierstrass semigroups. Let $P\in X$. In the case of $1$-Weierstrass points, 
the set 
$$
G(P):= \{\epsilon_i(P)+1: 0\le i\le g-1\}
$$
is the complement (or the gaps) of a  
semigroup $H(P)$, the so called Weierstrass semigroup at $P$. $H(P)$ 
looks like
$$
H(P)=\{0<m_1(P)<\ldots<m_g(P)=2g<2g+1<\ldots \},
$$
and it is satisfied the following property (\cite{B}, \cite[Thm. 
1.1(ii)]{Oliv}). Let $\ell_i(P):= \epsilon_i(P)+1$. Then
\begin{equation}\label{oliv}
\ell_i(P)\le 2i-2,\qquad \mbox{for}\qquad i=2,\ldots,g-1,\qquad \ell_g(P)\le 
2g-1,
\end{equation}
provided $m_1(P)\ge 3$. Then if $X$ is classical and if we have equality in 
(\ref{weight}), $w(P)$ can be computed by the formula  
\begin{equation}\label{sum}
w(P)=\frac{3g^2+g}{2}-\sum_{m\in H(P), m\le 2g}\ m.
\end{equation}
Now in case of $\gamma$-hyperelliptic curves, $P$  a ramified point of 
$\pi$ and $p>2$, $H(P)$ fulfil the following 
properties (\cite[Lemma 3.4]{To1}): 
\begin{enumerate}
\item[(A)] $\gamma=\# \{\ell\in G(P): \ell\ {\rm even}\}$.
\item[(B)] $H(\pi(P))= \{\frac{h}{2}: h \in H(P), h\ {\rm even}\}$.
\end{enumerate}
(Note that property (B) implies $h\le 2\gamma+2$ for $h\in H(P)$, $h$ even. 
In particular if $X$ is classical and $g\ge 2\gamma+2$, then each ramified 
point of $\pi$ is a Weierstrass point of $X$.)
\medskip

\noindent {\it Proof of (i).} Follows from property (A) above 
and \cite[Lemma 3.1.2(ii)]{To2}.
\smallskip

\noindent {\it Proof of (ii).} Let $P\in X$ be as in (II). Since $p=0$ or 
$p>2g-2>2\gamma-2$, then $\tilde X$ is also a classical curve.  Thus  
from properties (A) and (B) we have that all the even positive non-gaps 
of $H$ belong to the following set:
$$
\{2\gamma+2i: i\in \{1,\ldots,g-\gamma\}\}.
$$
Hence
\begin{equation*}
\sum_{h\in H(P), h\ {\rm even},\ h\le 2g} h = g^2+g-\gamma^2-\gamma.\tag{$*$}
\end{equation*}
Let denote by $u_\gamma<\ldots<u_1$ the $\gamma$ odd non-gaps at $P$ in 
$[1,2g-1]$. According to (\ref{sum}) and ($*$), an upper bound for $w(P)$ 
corresponds to a lower bound for $\sum_{i=1}^{\gamma} u_i$. By  
\cite[Lemma 2.1]{To1}, $u_\gamma\ge 2g-4\gamma+1$. If $u_\gamma\ge 
2g-2\gamma-1$, any odd number in $[2g-2\gamma+1,2g-1]$ could be an odd 
non-gap at $P$. Hence in this case we have
$$
\sum_{i=1}^{\gamma} u_i\ge 2\gamma g -\gamma^2 -2\gamma.
$$
If $u_\gamma\le 2g-2\gamma-3$ (then $\gamma\ge 2$), it is easy to see 
that the minimum for $\sum_{i=1}^{\gamma} u_i$ is reached for the sequence 
$2g-(2i+5), 2g-1$, $i=1,\ldots, \gamma-1$. Hence in this case we have 
$$
\sum_{i=1}^{\gamma}u_\gamma\ge 2\gamma g -\gamma^2 -4\gamma+4.
$$
Then since $\gamma>1$ from $(*)$, the last inequality and (\ref{sum}) we 
obtain  (ii).
\smallskip

\noindent {\it Proof of (iii).} Let $P\in X$ and suppose that $P$ is not 
a ramified point of $\pi$. 
\begin{claim*}
Let $h$ be a non-gap at $P$. Then $h\ge g-2\gamma+1$.
\end{claim*}
\begin{remark*} 
Let $f\in k(X)$ and denote by $O(f)$ the degree of $f$. Then 
$O(f)$ is 
even provided $O(f)<g+1-2\gamma$ and $g\ge 4\gamma+2$. For $p=0$  this is 
a result due to Farkas \cite[Thm.2(iii)]{F} (see also 
\cite[Thm.V.1.9]{F-K}, Accola \cite[Lemma 
4]{A}) and in general is due to Stichtenoth \cite[Satz 2]{Sti}. The claim 
follows from this result but with an extra hypothesis on $g$. We will 
see that in the case that $f$ has just one pole one can avoid such a 
hypothesis. The claim is a particular case of \cite[Corollary 3.3(ii)]{To1}, 
and for the sake of completeness we state a proof of it.
\end{remark*}
\begin{proof} {\it (Claim.)} 
Suppose that $h<g-2\gamma+1$. Consider $K':=k(\tilde X).k(f)$, with ${\rm 
div}_\infty (f)=hP$. Then by Castelnuovo's inequality concerning 
subfields of $k(X)$ (see \cite{C}, \cite{Sti1}) we must have 
$K'=k(\tilde X)$. Thus there exists $\tilde f\in K(\tilde X)$ such that 
$f=\tilde f\circ \pi$ and we would have that $P$ is a totally ramified 
point of $\pi$, a contradiction.
\end{proof}

Let $g\ge 2\gamma$. We have
\begin{equation}\label{w}
w(P)=\sum_{i=g-2\gamma+1}^{g}(\ell_i-i),
\end{equation}
and then we consider two cases:
\begin{enumerate}
\item[(a)] There exists $\ell \in G(P)\cap [g-2\gamma+1,g]$.
\item[(b)] $[g-2\gamma+1,g]\cap \mathbb N \subseteq H(P)$.
\end{enumerate}
In the first case we have $\ell_{g-2\gamma+1}=g-2\gamma+j$ with 
$j\in\{1,\ldots,\gamma\}$. Then from (\ref{w}) and (\ref{oliv}) we get
$$
w(P)\le 
(2\gamma-1)+\sum_{i=g-2\gamma+2}^{g-1}(i-2)+(g-1)=
(2\gamma-1)g-(\gamma-1)(2\gamma+1).
$$
In the second case, due to the semigroup property of $H(P)$, $w(P)$ 
reachs its maximum whenever 
$G(P)=\{1,\ldots,g-2\gamma,2g-6\gamma+2,\ldots, 2g-4\gamma+1\}$. Then 
from (\ref{w}) we find
$$
w(P)\le 2\gamma g - 2\gamma(4\gamma-1).
$$
This finish the proof of the lemma.


\begin{thebibliography}{Fet 9}

\bibitem[A]{A} Accola, R.D.M.: {\it Strongly branched coverings of closed 
Riemann surfaces,} Proc. Amer. Math. Soc. {\bf 26} (1970), 315--322.

\bibitem[B]{B} Buchweitz, R.O.:``\"Uber  deformationem monomialer 
kurvensingularit\"aten und 
Weierstrasspunkte auf Riemannschen fl\"achen", Thesis, Hannover 
1976.

\bibitem[C]{C} Castelnuovo, G.: {\it Sulle serie algebriche di gruppi di 
punti appartenenti ad una curva algebrica,} Rendiconti della Reale 
Accademia dei Lincei (5) {\bf 15} (1906), 337--344. Reprinted in Memoria 
scelta, Zanichelli, Bologna, 509--517, 1937.

\bibitem[F]{F} Farkas, H.M.: {\it Remarks on automorphisms of compact 
Riemann surfaces,} Ann. of Math. Stud. {\bf 78} (1974), 121--144.

\bibitem[F-K]{F-K} Farkas, H.M.; Kra, I.: ``Riemann surfaces", Grad. 
Texts in Math. {\bf 71} (second edition) Springer-Verlag 1992.

\bibitem[G]{G} Garcia, A.: {\it Weights of Weierstrass points in double 
covering of curves of genus one or two,} Manuscripta Math. {\bf 55} 
(1986), 419--432.

\bibitem[H-O]{H-O} Homma, M.; Ommori, S.: {\it On the weight of higher 
order  Weierstrass points,} Tsukuba J. Math. {\bf 8} (1984), 189--198.

\bibitem[K]{K} Kato, T.: {\it Non-hyperelliptic Weierstrass points of 
maximal weight,} Math. Ann. {\bf 239}, (1979), 141--147.

\bibitem[Lak]{Lak} Laksov, D.: {\it Weierstrass points on curves,} 
Ast\'erisque {\bf 87-88}, (Soci\'et\'e Math\'ematique de France, Paris, 
1981), 221--247.

\bibitem[Oliv]{Oliv} Oliveira, G.: {\it Weierstrass 
semigroups and the canonical 
ideal of non--trigonal curves,}  Manuscripta Math. {\bf 71} (1991), 
431--450.

\bibitem[Pf]{Pf} Pflaum, U.: {\it The canonical constellations of 
k-Weierstrass points,} Manuscripta Math. {\bf 59} (1987), 21--34.

\bibitem[Sti]{Sti} Stichtenoth, H.: {\it s-Erweiterungen algebraischer 
Funktionenk\"orper,} Arch. Math. {\bf 43} (1984), 27--31.

\bibitem[Sti1]{Sti1} Stichtenoth, H.: {\it Die Ungleichung von 
Castelnuovo,} J. Reine Angew. Math. {\bf 348} (1984), 197--202.

\bibitem[S-V]{S-V} St\"ohr, K.O.; Voloch, J.F.: {\it Weierstrass points 
and curves over finite fields,} Proc. London Math. Soc. (3), {\bf 52} 
(1986), 1--19.

\bibitem[To1]{To1} Torres, F.: {\it On certain $N$-sheeted coverings of 
curves and numerical semigroups which cannot be realized as Weierstrass 
semigroups,} Comm. Algebra {\bf 23} (11) (1995), 4211--4228.

\bibitem[To2]{To2} Torres, F.: {\it Remarks on numerical semigroups,} 
alg-geo e-print 9512012.


\end{thebibliography}
\end{document}